\documentclass[11pt,twoside,a4paper]{article}
\usepackage{amsfonts}
\usepackage{amssymb}
\usepackage{amsmath}

\begin{document}

\newcommand{\pa}{\partial}
\newcommand{\opa}{\overline\pa}
\newcommand{\ol}{\overline}

\numberwithin{equation}{section}

\newcommand\C{\mathbb{C}}  
\newcommand\R{\mathbb{R}}
\newcommand\Z{\mathbb{Z}}
\newcommand\N{\mathbb{N}}
\newcommand\PP{\mathbb{P}}

{\LARGE \centerline{On the normal bundle of Levi-flat real hypersurfaces}}
\vspace{0.8cm}

\centerline{\textsc {\large Judith Brinkschulte}\footnote{Universit\"at Leipzig, Mathematisches Institut, PF 100920, D-04009 Leipzig, Germany. 
E-mail: brinkschulte@math.uni-leipzig.de
 }}

\vspace{0.5cm}

\begin{abstract} Let $X$ be a connected complex manifold of dimension $\geq 3$ and $M$ be a smooth compact Levi-flat real hypersurface in $X$.  We show that the normal bundle to the Levi foliation does not admit a Hermitian metric with positive curvature along the leaves. This generalizes a result obtained by Brunella.
\end{abstract}

\vspace{0.5cm}

\section{Introduction}

A real hypersurface $M$ (of class at least $\mathcal{C}^2$) in a complex manifold is called Levi-flat if its Levi-form vanishes identically or, eqivalently, if it admits a foliation by complex hypersurfaces. Another equivalent formulation is that $M$  is locally pseudoconvex from both sides.\\

Given a Levi-flat real hypersurface $M$ in a complex manifold $X$ of dimension $n$, we call $N^{1,0}_M = (T^{1,0}_X)_{\mid M} / T^{1,0}M$ the holomorphic normal bundle of $M$. The restriction of $N^{1,0}_M$ to each $(n-1)$-dimensional complex submanifold of $M$ has a structure of a holomorphic line bundle induced from that of $T^{1,0}_X$.\\

In this paper we prove the following\\

\newtheorem{main}{Theorem}[section]
\begin{main}   \label{main}   \ \\
Let $X$ be a complex manifold of dimension $n \geq 3$. Then there does not exist a smooth compact Levi-flat real hypersurface $M$ in $X$ such that the normal bundle to the Levi foliation  admits a  Hermitian metric with positive curvature along the leaves.
\end{main}

Classical nontrivial examples of Levi-flat hypersurfaces were described by Grauert as tubular neighborhoods of the zero section of a generically chosen line bundle over a non-rational Riemann surface \cite{G}. In these examples, the Levi-flat hypersurfaces arise as the boundary of a pseudoconvex domain admitting only constant holomorphic functions. On the other hand, there are also examples of compact Levi-flat real hypersurfaces bounding Stein domains. For example, the product of an annulus and the punctured plane is bilomorphically equivalent to a domain in $\PP^1\times \lbrace\C/ (\mathbb{Z} + i\mathbb{Z})\rbrace$ with Levi-flat boundary \cite{O1}.  Further examples of complex surfaces that can be cut into two Stein domains along smooth Levi-flat real hypersurface can be found in \cite{N}. From \cite{Ad} one even obtains examples of Levi-flat hypersurfaces in  complex surfaces having hyperconvex complement.\\

These  examples above show that Levi-flat hypersurfaces can be of quite different nature and therefore explain a certain interest in the classification of compact Levi-flat real hypersurfaces. Let us also mention that some of these constructions can be extended to higher dimensions.\\

On the other hand, the study of Levi-flat real hypersurfaces is related to basic questions in dynamical systems and foliation theory: Levi-flats arise as stable sets of holomorphic foliations, and a real-analytic Levi-flat real hypersurface extends to a holomorphic foliation leaving $M$ invariant. Relating to this, a famous open problem is whether or not $\C\PP^2$ contains a smooth Levi-flat real hypersurface. This problem arose as part of a conjecture that, for any codimension one holomorphic foliation on $\C\PP^2$ (with singularities), any leaf accumulates to a singular point of the foliation \cite{CLS}. This problem is still open. It is only known that if $\C\PP^2$ admits a smooth Levi-flat real hypersurface, then it has to satisfy a restrictive curvature condition \cite{AB}.\\

For $n\geq 3$, however, it is known that there does not exist any smooth real Levi-flat hypersurface $M$  in $\C\PP^n$. This was first proved by LinsNeto in $\cite{LN}$ for real-analytic $M$ and by Siu in $\cite{S}$ for $\mathcal{C}^{12}$-smooth $M$. For further improvements concerning the regularity, we refer the reader to $\cite{IM}$ and $\cite{CS}$. \\

The proofs of the above-mentioned results essentially exploited the positivity of $T^{1,0}\C\PP^n$.
Brunella's main observation \cite{Br} was that the positivity of the normal bundle itself is enough to ensure that the complement of $M$ is pseudoconvex. If $X= \C\PP^n$, or if $X$ admits a hermitian metric of positive curvature, then the normal bundle $N^{1,0}_M$ is automatically positive (it is a quotient of $T^{1,0}X$, and therefore more positive than $T^{1,0}X$).\\

This led Brunella  to prove that if $X$ is a compact K\"ahler manifold with $\dim X \geq 3$, and if $M$ is a smooth Levi-flat real hypersurface such that there exists a holomorphic foliation on a neighborhood of $M$ leaving $M$ invariant, then the normal bundle of this foliation does not admit any fiber metric with positive curvature. \\

Ohsawa generalized this in \cite{O7} to a nonexistence result for smooth Levi-flat real hypersurfaces admitting a fiber metric whose curvature form is semipositive of rank $\geq 2$ along the leaves of $M$ (in any compact K\"ahler manifold).\\

Our Theorem \ref{main} is a generalized version of Brunella's result in the sense that we are able to drop the compact K\"ahler assumption on the ambient $X$. This was conjectured in \cite[Conjecture\ 5.1]{O8}. \\

The following example from \cite[Example\ 4.2]{Br} and \cite[Theorem\ 5.1]{O8} shows that Theorem \ref{main} cannot hold for $n=2$, even for $X$ compact K\"ahler:\\

Let $\Sigma$ be a compact Riemann surface of genus $g\geq 2$. Let $\mathbb{D}$ be the open unit disc, and let $\Gamma$ be a discrete subgroup of $\mathrm{Aut}\mathbb{D}\subset\mathrm{Aut}\C\PP^1$ such that $\Sigma\simeq \mathbb{D}/\Gamma$. Then $\Gamma$ also acts on $\mathbb{D}\times\C\PP^1$ by
$$(z,w)\mapsto (\gamma(z),\gamma(w)),\ \gamma\in\Gamma.$$
The quotient $X= (\mathbb{D}\times\C\PP^1)/\Gamma$ is a compact complex surface, ruled over $\Sigma$ (and hence projective). From the horizontal foliation on $\mathbb{D}\times\C\PP^1$, we get a holomorphic foliation on $X$, leaving invariant a real analytic Levi-flat hypersurface $M$ induced from the $\Gamma$-invariant $\mathbb{D}\times \mathbb{S}^1$. The Bergman metric induces a metric with positive curvature on the normal bundle of $M$ (see \cite{O8} for more details).\\

\vspace{1cm}

{\bf Acknowledgements.} I wish to express my thanks to Stefan Nemirovski for his contribution to this paper: the construction of the K\"ahler metric in the proof of Proposition \ref{kähler} was essentially his idea. I would also like to thank Masanori Adachi and Takeo Ohsawa not only for their great interest, but also for  many helpful remarks improving the paper.\\
The research on this project was supported by Deutsche Forschungsgemeinschaft (DFG, German Research Foundation, grant BR 3363/2-2).\\

\section{Sketch of the proof}

Let us begin by recalling the essential steps of Brunella's proof in \cite{Br}: Assume that $X$ is a connected compact K\"ahler manifold of dimension $n\geq 3$, and let $M$ be a smooth Levi-flat real hypersurface such that there exists a holomorphic foliation on a neighborhood of $M$ leaving $M$ invariant.
Under the assumption that the normal bundle of this foliation admits a fiber metric with positive curvature,  Brunella shows that $X\setminus M$ is strongly pseudoconvex. Then the argument is as follows: Since the normal bundle of the foliation is topologically trivial, its curvature form $\theta$ is $d$-exact on a tubular neighborhood $U$ of $M$. Thus $\theta = d\beta$ on $U$, where the primitive $\beta = \beta^{1,0} + \beta^{0,1}$ can be chosen of real type ($\ol\beta^{1,0} = \beta^{0,1}$) and one has $\opa\beta^{0,1}=0$. Since $\dim X \geq 3$, the vanishing theorem of Gauert and Riemenschneider combined with Serre's duality implies that the $\opa$-cohomology with compact support $H_c^{0,2}(X\setminus M)$ is zero. This means that one can extend $\beta^{0,1}$ $\opa$-closed to $X$. Hodge symmetry on the K\"ahler manifold $X$ means $H^{0,1}(X) \simeq \overline{H^{1,0}(X)}$. Hence $\beta^{0,1} = \eta + \opa\alpha$, with $\pa\eta =0$. But then $\pa\beta^{0,1} = \pa\opa\alpha$. Therefore, setting $\phi = i(\ol\alpha - \alpha)$, one obtains $\theta = i\pa\opa\phi$. The existence of a potential for the positive curvature form is, however, a contradiction to the maximum principle on the leaves of the foliation.\\

Our proof follows this general idea. We assume by contradiction that there exists a smooth compact Levi-flat real hypersurface $M$ in $X$ such that the normal bundle to the Levi foliation  admits a Hermitian metric with positive curvature along the leaves.
 However, since our $M$ is not embedded in a compact K\"ahler manifold, we have to make several important modifications. Since $M$ has a tubular neighborhood which is pseudoconcave (of dimension $\geq 3$),  this tubular neighborhood can be compactified to a compact manifold $\tilde X$. Then $\tilde X\setminus M$ is a strongly pseudoconvex manifold, and we can even arrange that it carries a complete K\"ahler metric (section  4 and 5). By means of $L^2$-estimates on $\tilde X \setminus M$, we will then extend the normal bundle to $M$ to a holomorphic line bundle over $\tilde X$ (section 6). We also show that $CR$ sections of high tensor powers of the normal bundle extend to holomorphic sections over $\tilde X$ (section 7), again by means of solving some Cauchy-problem for the $\opa$-equation using $L^2$-estimates. This permits us to find sufficiently many sections that provide a holomorphic embedding of a tubular neighborhood of $M$ into a compact K\"ahler manifold (section 8). This proves the nonexistence of such $M$ as before.\\

\section{Preliminaries}

Let $Y$ be a complex manifold of dimension $n$ endowed with a Hermitian metric $\omega$, and let $E$ be a holomorphic vector bundle on $Y$ with a Hermitian metric $h$. For integers $0\leq p,q\leq n$, we use the following notations:\\
$\mathcal{C}^{p,q}_c(Y,E)$ denotes the space of smooth, compactly supported $E$-valued $(p,q)$-forms on $Y$.\\
$L^2_{p,q}(Y,E,\omega,h)$ denotes the Hilbert space obtained by completing $\mathcal{C}^{p,q}_c(Y,E)$ with respect to the $L^2$-norm $\Vert \cdot\Vert_{\omega,h}$ induced by $\omega$ and $h$.\\

If $\mathrm{rk}E =1$, and the metric on $E$ is given by $h= e^{-\varphi}$, we write $L^2_{p,q}(Y,E,\omega,\varphi)$ instead of $L^2_{p,q}(Y,E,\omega,h)$.\\

As usual, the differential operator $\opa$ is extended as a densely defined closed linear operator on $L^2_{p,q}(Y,E,\omega,h)$, whose domain of definition is 
$$\mathrm{Dom}\opa = \lbrace f\in L^2_{p,q}(Y,E,\omega,h)\mid \opa f \in L^2_{p,q+1}(Y,E,\omega,h)\rbrace,$$
where $\opa f$ is computed in the sense of distributions. The Hilbert space adjoint of $\opa$ will be denoted by $\opa^\ast \ (= \opa^\ast_{\omega,h})$.\\

We also define the space of harmonic forms,
$$\mathcal{H}^{p,q}(Y,E,\omega,h)= L^2_{p,q}(Y,E,\omega,h)\cap\mathrm{Ker}\opa\cap\mathrm{Ker}\opa^\ast_{\omega,h},$$
 and the $L^2$-Dolbeault cohomology groups of $Y$,
$$H^{p,q}_{L^2}(Y,E,\omega,h) = L^2_{p,q}(Y,E,\omega,h)\cap\mathrm{Ker}\opa / L^2_{p,q}(Y,E,\omega,h)\cap\mathrm{Im}\opa.$$

Whenever we feel that it is clear from the context, we will omit the dependency of the $L^2$-spaces, norms, operators etc. on the hermitian metric $h$ of the vector bundles under considerations. \\

For the reader's convenience, we also recall the well-known Bochner-Kodaira-Nakano inequality, which is the starting point for all $L^2$-estimates for $\opa$.\\

If $\omega$ is a K\"ahler metric, then for every $u\in \mathcal{C}^{p,q}_c(Y,E)$ we have the following a priori estimate (see $\cite[\mathrm{Lemme}\ 4.4]{D}$):
\begin{equation}   \label{BKN}
 \Vert \opa u\Vert^2_{\omega,h} + \Vert \opa^\ast_{\omega,h} u\Vert^2 \geq \ \ll \lbrack i \Theta_h(E),\Lambda_\omega\rbrack u,u\gg_{\omega,h}
\end{equation}

Here $i\Theta_h(E)$ is the curvature of the bundle $E$, and $\Lambda_\omega$ is the adjoint of multiplication by $\omega$.  It is important to note that if the metric $\omega$ is {\it complete}, then the inequality (\ref{BKN}) extends to all forms $u\in L^2_{p,q}(Y,E,\omega,h)\cap\mathrm{Dom}\opa\cap\mathrm{Dom}\opa^\ast_{\omega,h}$. For metrics that are not K\"ahler, there is an additional curvature term (see \cite{O8}).\\

Moreover, if $E$ is a line bundle, and if $\lambda_1\leq \ldots\leq\lambda_n$ are the eigenvalues of $i\Theta(E)$ with respect to $\omega$, then we have
\begin{equation}  \label{curvature}
\langle \lbrack i\Theta_h(E),\Lambda_\omega\rbrack u,u\rangle_{\omega,h} \geq(\lambda_1 +\ldots +\lambda_q -\lambda_{p+1} - \ldots - \lambda_n)\vert u\vert^2_{\omega,h}
\end{equation}
if $u$ is of bidegree $(p,q)$ (see $\cite[(13.6)]{D2}$).\\

In section 7, we shall also use the following variant of the $\opa$-operator: by $\opa_c$ we denote the strong minimal realization of $\opa$ on $L^2_{p,q}(Y,E,\omega,h)$. This means that $u\in\mathrm{Dom}\opa_c\subset L^2_{p,q}(Y,E,\omega,h)$ if there exists $f\in L^2_{p,q+1}(Y,E,\omega,h)$ and a sequence $(u_\nu)_{\nu\in\N} \subset\mathcal{D}^{p,q}(Y,E)$ such that $u_\nu\longrightarrow u$ and $\opa u_\nu\longrightarrow f = \opa_c u$ in $L^2_{p,q+1}(Y,E,\omega,h)$.\\

The Hilbert space adjoint of $\opa_c$ will be denoted by $\vartheta$; it is the weak maximal realization of the formal adjoint of $\opa$ on $L^2_{p,q}(Y,E,\omega,h)$.

\section{Convexity properties}

Let $M$ be a smooth Levi-flat real hypersurface in a Hermitian manifold $X$. By considering a double covering, we may assume that $M$ is orientable and that the complement of $M$ in $X$ divides $X$ into two connected components (shrinking $X$ if necessary), see also \cite{Br}. So we may assume that $X$ is sufficiently small so that there exists a smooth real-valued function $\rho$ on $X$ such that 
$$M = \lbrace z\in X\mid \rho(z)=0\rbrace$$
 and $d\rho \not=0$ on $M$.
We further fix a Hermitian metric  $\omega_o$  on $X$.\\

 We now assume by contradiction that the normal bundle of the Levi foliation admits a Hermitian metric of positive curvature along the leaves. As in \cite{Br}, this implies that the complement of $M$ is strongly pseudoconvex. The following proposition was proved in \cite{O7}, we include the proof for the sake of completeness.\\

\newtheorem{convexity}[main]{Proposition}
\begin{convexity}     \label{convexity}  \ \\
Let $M$ be a compact smooth Levi-flat real hypersurface in a Hermitian manifold $X$ of dimension $\geq 2$, such that the normal bundle $N_M^{1,0}$ of the Levi foliation admits a smooth Hermitian metric of leafwise positive curvature. Then, after possibly shrinking $X$, there exists  a $\mathcal{C}^2$-smooth nonnegative  function $v$  on $X$, smooth on $X\setminus M$,  and a positive constant $c > 0$ such that
\begin{equation}  \label{log}
i\pa\opa(-\log v) \geq c\omega_o\qquad\mathrm{on}\ X\setminus M.
\end{equation}
Moreover, we have that $v = O(\rho^2)$
\end{convexity}

{\it Proof.} 
As in \cite{Br} one can find a finite family of local coordinate neighborhoods $U_\gamma$ in $X$ covering $M$ such that $U_\gamma \cap M = \lbrace z\in U_\gamma\mid \mathrm{Im} f_\gamma =0\rbrace$, where $\opa f_\gamma$ vanishes to infinite order on $M$, and that $T^{1,0}_M = \mathrm{Ker} d f_\gamma$ on $M\cap U_\gamma$. \\

Let $\varpi = \lbrace \varpi_\gamma\rbrace$ be a $1$-form on $M$ defining its Levi-foliation. We may assume that $\varpi_\gamma$ is defined on $U_\gamma$. Let $h = \lbrace h_\gamma\rbrace$ be the fibre metric of $N^{1,0}_M$ such that 
\begin{equation}  \label{1}
h_\delta = \vert \varpi_\gamma / \varpi_\delta\vert^2 h_\gamma
\end{equation}
on $U_\gamma\cap U_\delta\cap M$. By assumption on the curvature of $N^{1,0}_M$, we may assume that $-\log h_\gamma$ is strictly plurisubharmonic on the leaves of $M$. \\

We have $\varpi_\gamma = e_\gamma d f_\gamma$ for some smooth function $e_\gamma$ which is nowhere vanishing on $U_\gamma$ and holomorphic along the leaves of $M$. From (\ref{1}) it follows that we have
\begin{equation*}
h_\gamma \vert e_\gamma\vert^2 (\mathrm{Im}f_\gamma)^2 - h_\delta \vert e_\delta\vert^2 (\mathrm{Im}f_\delta)^2 = O(\rho^3) \ \mathrm{on}\ U_\gamma\cap U_\delta
\end{equation*}

Therefore, invoking Whitney's extension theorem, there exists a $\mathcal{C}^2$ function $v$ defined in a tubular neighborhood of $M$, smooth away from $M$, such that $v = h_\gamma \vert e_\gamma\vert^2 (\mathrm{Im}f_\gamma)^2 + O(\rho^3)$ on $U_\gamma$. \\

To see that $-\log v$ is strictly plurisubharmonic in a tubular neighborhood of $M$, it suffices to estimate the Levi-form of $-\log(h_\gamma \vert e_\gamma\vert^2 (\mathrm{Im}f_\gamma)^2) $. Indeed, let $V\in T^{1,0}X$ be a unitary vector that we decompose orthogonally into $V = V_t + V_n$, with $V_t\in\mathrm{Ker}\pa\rho$.  Then the strict plurisubharmonicity of $-\log h_\gamma$ and the holomorphicity of $e_\gamma$ along the leaves of $M$ imply that there exists $c > 0$ such that
$$-i\pa\opa \log(h_\gamma \vert e_\gamma\vert^2) (V_t,\ol V_t) \geq (c-\epsilon) \omega_o( V_t,\ol V_t),$$
where $\epsilon$ can be made as small as we wish by shrinking $X$. On the other hand, since $h_\gamma$ and $\vert e_\gamma\vert^2$ do not vanish, 
$$-i\pa\opa \log(h_\gamma \vert e_\gamma\vert^2) (V_n,\ol V_n) \geq -C \omega_o(V_n,\ol V_n).$$
The mixed terms in $(V_t,\ol V_n)$ can be handled as follows:
$$-i\pa\opa \log(h_\gamma \vert e_\gamma\vert^2) (V_t,\ol V_n) \geq -\epsilon \omega_o( V_t,\ol V_t) - \frac{C}{\epsilon} \omega_o(V_n,\ol V_n).$$

Moreover, since $\opa f_\gamma$ vanishes to infinite order on $M$, we have 
\begin{eqnarray*}
-i\pa\opa\log ((\mathrm{Im}f_\gamma)^2) (V,\ol V) & = &(-2 \frac{i\pa\opa\mathrm{Im} f_\gamma}{\mathrm{Im}f_\gamma} + 2i \frac{\pa\mathrm{Im}f_\gamma\wedge\opa\mathrm{Im}f_\gamma}{(\mathrm{Im}f_\gamma)^2}) (V,\ol V)\\
& \geq & -\epsilon \omega_o(V,\ol V) + 2i \frac{\pa\mathrm{Im}f_\gamma\wedge\opa\mathrm{Im}f_\gamma}{(\mathrm{Im}f_\gamma)^2}(V_n,\ol V_n)
\end{eqnarray*}
Again, $\epsilon$ can be made as small as we wish by shrinking $X$. Combining the above estimates permits to conclude by taking $\epsilon$ sufficiently small. \hfill$\square$\\

{\it Remark.} In \cite[Proposition\ 3.3]{Ad2}, a converse statement is proved: If there exists a boundary distance function of $X\setminus M$ with positive Diederich-Fornaess exponent, then the normal bundle $N^{1,0}_M$ is positive along the leaves.\\

\section{A first compactification}

For sufficiently large $\alpha \in\R^+$, Proposition \ref{convexity} implies that the  set\\
 $\lbrace z\in X\mid -\log v(z) > \alpha\rbrace $ is a pseudoconcave manifold (of complex dimension $\geq 3$). By a theorem of Rossi \cite{R} and Andreotti-Siu \cite{AS}
it can be compactified. Hence we may assume that $M$ is embedded as a real hypersurface in a compact complex manifold $X^\prime$ of dimension $n$, and $X^\prime\setminus M$ is a strongly pseudoconvex manifold (or a $1$-convex manifold, using a different terminology): $X^\prime\setminus M$ admits a smooth exhaustion function (given by Proposition \ref{convexity}), which is strictly plurisubharmonic outside a compact subset.\\

Before continuing with the proof, we will make some standard modifications of $X^\prime$ in order to faciliate the following arguments.\\

By \cite{G1} there exists a compact analytic subset $A\subset X^\prime\setminus M$ and a proper holomorphic map $\pi^\prime$ from $X^\prime\setminus M$ onto a Stein space $S$ such that $\pi^\prime$ is a biholomorphic mapping from $X^\prime\setminus(M\cup A) $ to $S\setminus \pi^\prime (A)$. \\

By Hironaka's method, there is a complex manifold $\tilde S$ obtained from $S$ by blowing up $S$ along smooth centers, several times if necessary, such that the induced bimeromorphic map $\tilde\pi: \tilde S\rightarrow S$ is holomorphic. Moreover, following $\cite{G1}$, it is possible to choose $\tilde S$ such that 
\begin{enumerate}
\item[$\bullet$] $(\pi^\prime)^{-1}\circ\tilde\pi$ is biholomorphic on $\tilde S\setminus \tilde\pi^{-1}(\pi^\prime(A))$,
\item[$\bullet$] $\tilde\pi^{-1}(\pi^\prime(A))$ is a divisor with normal crossings whose irreducible components $\lbrace \tilde A_j\rbrace_{j=1}^{\nu}$ are nonsingular, and
\item[$\bullet$] there exists positive integers $p_1,\ldots,p_\nu$ such that the line bundle $\mathcal{O}(\tilde A)$ induced by the divisor  $\tilde A = \sum_{j=1}^{\nu} p_j \tilde A_j$  is negative. 
\end{enumerate}

This modification permits us to prove the following proposition.\\

\newtheorem{compactify}[main]{Proposition}
\begin{compactify}     \label{compactify}  \ \\
$M$ is a real hypersurface in a compact complex manifold $\tilde X$ of dimension $n$ such that
\begin{enumerate}
\item[(i)] $\tilde X\setminus M$ is strongly pseudoconvex, and moreover there exists a smooth exhaustion function $\varphi$ on $\tilde X\setminus M$, plurisubharmonic on $\tilde X\setminus M$, strictly plurisubharmonic outside a compact of $\tilde X\setminus M$, such that $\varphi \sim -2\log \vert\rho\vert$ outside a compact of $K$ of $\tilde X\setminus M$.
\item[(ii)] $\tilde X\setminus M$ admits a complete K\"ahler metric $\tilde\omega$ such that $\tilde\omega = i\pa\opa \varphi$ outside a compact of $\tilde X\setminus M$.
\item[(iii)] There exists a line bundle $L$ over $\tilde X$, defined by a divisor $\tilde A = \sum_{j=1}^\nu p_j \tilde{A_j}$, with compact support in $\tilde X\setminus M$, such that $L$ is negative on an open neighborhood of $K$ and
holomorphically trivial outside a compact of $\tilde X\setminus M$.

\end{enumerate}
\end{compactify}

{\it Remark.} In the following sections, we implicitely assume that $L$ is endowed with a flat metric outside of $K$ (where $L$ is holomorphically trivial).\\

{\it Proof.}
Gluing a pseudoconcave tubular neighborhood of $M$ in $X$ to a suitable relatively compact domain with strictly pseudoconvex boundary in $\tilde S$ we obtain a compact complex manifold $\tilde X$ of dimension $n$, containing the Levi-flat real hypersurface $M$, that has all the required properties. Indeed, (i) and (iii) follow easily from the discussion preceeding the proposition. \\

To see that $\tilde X\setminus M$ is complete K\"ahler, we consider the K\"ahler metric $\omega_{\tilde S} = i\Theta (L^\ast) = i\Theta(\mathcal{O}(-\tilde A))$ on $\tilde S$. On $\tilde S\setminus \tilde A$, the line bundle $\mathcal{O}(-\tilde A)$ is holomorphically trivial, hence there exists a smooth function $\psi$ such that $i\pa\opa\psi = i\mathcal{O}(-\tilde A)$ on $\tilde S\setminus \tilde A$. 
Now we choose a smooth cut-off function $\chi$ such that $\chi \equiv 1$ on a sufficiently large compact of $\tilde S$ containing $\tilde A$ and such that the support of $\chi$ is compact in $\tilde X\setminus M$. 
Then, for $\varepsilon > 0$ small enough,
\begin{equation} \nonumber
\tilde \omega = \left\{ 
\begin{aligned}
\varepsilon i\pa\opa(\chi\psi) + i\pa\opa \varphi \\
\varepsilon i\Theta(\mathcal{O}(-\tilde A)) + i\pa\opa\varphi  
\end{aligned}   \right.
\end{equation}
defines a K\"ahler metric on $\tilde X\setminus M$. Moreover, $\tilde \omega$ is complete on $\tilde X\setminus M$. 
Indeed, it  follows from (\ref{log}) as  in $\cite{OS1}$ that there exists $0< \eta\leq 1$ such that $i\pa\opa (-v^\eta) \gtrsim v^\eta \omega_o$. This implies that $i\pa\opa (-\log v) \gtrsim i\eta \pa\log v \wedge\opa\log v$, showing that $\tilde\omega$ is complete.\hfill$\square$\\

\section{Holomorphic extension of the normal bundle}

The aim of this section is to prove that the holomorphic normal bundle of $M$ extends to a holomorphic line bundle over the compact manifold $\tilde X$. The main ingredient needed for the extension is the following $L^2$-vanishing result:\\

\newtheorem{vanish1}{Proposition}[section]
\begin{vanish1}     \label{vanish1}  \ \\
For every $N\in \N$ the following holds: Assume $0\leq q\leq n-1$, and let $f\in L^2_{0,q}(\tilde X\setminus M,\tilde\omega, -N\varphi)\cap\mathrm{Ker}\opa$. Then there exists a $(0,q-1)$-form $g\in L^2_{0,q-1}(\tilde X\setminus M,\tilde\omega, -N\varphi)$ satisfying $\opa g = f$.
\end{vanish1}

{\it Proof.} The proof is similar to the one of Theorem 2.1 in \cite{O7}. The metric $\tilde\omega$ is K\"ahler, and, since $\varphi$ is plurisubharmonic on $\tilde X\setminus M$, we have $-i\pa\opa\varphi\leq 0$ on $\tilde X\setminus M$ and $-i\pa\opa\varphi = \tilde\omega$ outside a compact $K$ of $\tilde X\setminus M$. From (\ref{BKN}) and (\ref{curvature}) we then obtain 
\begin{equation}  \label{outsideK}
N\int_{\tilde X\setminus (M\cup K)} \vert u\vert^2_{\tilde\omega} e^{N\varphi} dV_{\tilde\omega} \leq \Vert \opa u\Vert^2_{\tilde\omega, -N\varphi} + \Vert \opa^\ast_{\tilde\omega,-N\varphi} u\Vert^2_{\tilde\omega, -N\varphi} 
\end{equation} 
for every $u\in L^2_{0,q}(\tilde X\setminus M,\tilde \omega, -N\varphi) \cap \mathrm{Dom}\opa\cap\mathrm{Dom}\opa^\ast_{\tilde\omega,-N\varphi}$, $0\leq q\leq n-1$, $N\in\N$. \\

It is well-known that (\ref{outsideK}) implies that $\mathrm{Im}\opa$ is closed in $L^2_{0,q}(\tilde X\setminus M,\tilde \omega, -N\varphi)$ and that $H^{0,q}_{L^2}(\tilde X\setminus M, \tilde \omega, -N\varphi)$ is finite dimensional. This entails 
$$H^{0,q}_{L^2}(\tilde X\setminus M, \tilde \omega, -N\varphi)\simeq \mathcal{H}^{0,q} (\tilde X\setminus M, \tilde \omega, -N\varphi).$$
By (\ref{outsideK}) it then follows that every element of $\mathcal{H}^{0,q} (\tilde X\setminus M, \tilde \omega, -N\varphi)$ is zero outside of $K$, so that it vanishes identically by Aronszajn's unique continuation theorem for elliptic operators. Hence $H^{0,q}_{L^2}(\tilde X\setminus M, \tilde \omega, -N\varphi)= \lbrace 0\rbrace$ .\hfill$\square$\\

Since close to $M$, the weight  $e^{N\varphi}$ (up to a bounded function) equals $\rho^{-2N}$ (where $\rho$ is a defining function for $M$, see section 4), Proposition \ref{vanish1} enables us to extend $CR$ objects on $M$ to holomorphic objects on $\tilde X$ by solving $\opa$-equations with zero Cauchy data along $M$. In particular we can prove the following

\newtheorem{extensionbdle}[vanish1]{Proposition}
\begin{extensionbdle}     \label{extensionbdle}  \ \\
There exists a holomorphic line bundle $\tilde N\longrightarrow \tilde X$ such that $\tilde N_{\mid M} = N^{1,0}M$ ($\tilde N$ extends the $CR$ line bundle $N^{1,0}_M$). 
\end{extensionbdle}

{\it Proof.} The $CR$ line bundle $N^{1,0}_M$ is topologically trivial over $M$. Therefore it is in the image of the exponential map 
$$\exp :\  H^1(M,\mathcal{O}_M) \longrightarrow H^1(M,\mathcal{O}^\ast_M),$$
where $\mathcal{O}_M$, resp $\mathcal{O}_M^\ast$ denotes the sheaf of germs of smooth $CR$ functions on $M$, resp. nonvanishing $CR$ functions on $M$. Let us therefore choose $\xi\in H^1(M,\mathcal{O}_M)$, and identify it with a smooth $\opa_M$-closed $(0,1)$-form $g$ on $M$. Now $g$ admits a $\opa$-closed extension to $\tilde X$. Indeed, consider a smooth extension $\tilde g$ to a neighborhood of $M$ such that $\opa\tilde g$ vanishes to infinite order along $M$.
 Multiplying $\tilde g$ by a cut-off function whose support is contained in an arbitrary small tubular neighborhood of $M$, we can arrange that $\opa\tilde g$ vanishes outside a small tubular neighborhood of $M$. 
This means that the $(0,2)$-form $f = \opa\tilde g$ satisfies 
$$\int_{\tilde X\setminus M}\vert f\vert^2_{\tilde \omega} e^{N\varphi} dV_{\tilde\omega} < + \infty$$
for any $N\in\N$. Since $2\leq n-1$ by assumption on $X$, we may apply Proposition \ref{vanish1} and obtain
a smooth $(0,1)$-form $u$ satisfying $\opa u = \opa\tilde g$ on $\tilde X\setminus M$ and 
$$\int_{\tilde X\setminus M}\vert u\vert^2_{\tilde\omega} e^{N\varphi} dV_{\tilde\omega} < + \infty.$$
The solution that is minimal with respect to this $L^2$-norm moreover satisfies an elliptic equation. Using the regularity result obtained in $\cite[\mathrm{Theorem}\ 2.1]{B}$, we may therefore assume  that $u$ vanishes to sufficiently high order along $M$ (by taking $N$ sufficiently large). But then $\tilde g - u$ is $\opa$-closed on $\tilde X$ and coincides with $g$ on $M$. Therefore the holomorphic line bundle defined by $\tilde N=\exp(\tilde g - u)$ extends the $CR$ line bundle $N^{1,0}_M$. It is topologically trivial, since it is in the image of the exponential map above. \hfill$\square$\\

\section{Holomorphic extension of sections}

The key result of this section is Proposition \ref{extension}, the extension of $CR$ sections of the normal bundle
to holomorphic sections of $\tilde N$ over $\tilde X$. This will enable us to holomorphically embed a tubular neighborhood of $M$ into some complex projective space in the next section. The proof of this holomorphic extension property needs several steps; it is actually the technically most demanding part of this paper.\\

In order to extend $CR$ sections of a  line bundle over $M$ to holomorphic sections over $\tilde X$, the following vanishing result is very useful; it is in the same spirit as  Proposition \ref{vanish1}, but less precise. Remember that the holomorphic line bundle $L\longrightarrow \tilde X$ is given by Proposition \ref{compactify}.

\newtheorem{vanish2}{Proposition}[section]
\begin{vanish2}     \label{vanish2}  \ \\
Let $E\longrightarrow\tilde X$ be a holomorphic line bundle. Then there exist $k \in \N$  and $N\in\N$ such that the following holds: Assume $0\leq q\leq n-1$, and let $f\in L^2_{0,q}(\tilde X\setminus M, E\otimes L^k, \tilde\omega, -N\varphi)\cap\mathrm{Ker}\opa$. Then there exists a $(0,q-1)$-form $g\in L^2_{0,q-1}(\tilde X\setminus M, E\otimes L^k, \tilde\omega, -N\varphi)$ satisfying $\opa g = f$ and $\Vert g\Vert_{\tilde\omega, -N\varphi}\leq\Vert f\Vert_{\tilde\omega,-N\varphi}$.
\end{vanish2}

{\it Proof.} Given the holomorphic line bundle $E$, we first choose $k$ big enough such that $E\otimes L^k$ is negative on the compact where $\varphi$ is only weakly plurisubharmonic. Then we choose $N$ big enough such that $i\Theta(E\otimes L^k)-Ni\pa\opa\varphi \leq -\tilde\omega$ on $\tilde X$. 
We may then conclude from (\ref{BKN}) and (\ref{curvature}) as in the proof of Proposition \ref{vanish1}.
\hfill$\square$\\

In the next section, we want to holomorphically extend $CR$ sections of some high tensor power of the normal bundle of $M$. Since the normal bundle $N^{1,0}_M$ is positive, also $\tilde N$ can be equipped with a metric of  positive curvature near $M$.\\
On the other hand, we can multiply the metric of $\tilde N$ by $e^{N\varphi}$. This adds $-Ni\pa\opa\varphi$ to the curvature, so the curvature of $\tilde N$ can be made negative near $M$ by taking $N$ sufficiently large. This modification, however, would require the $CR$ sections that we wish to extend to be sufficiently regular. In \cite{Ad} it was shown that even $CR$ sections given an embedding into projective space are not necessarily $\mathcal{C}^\infty$-smooth if $\dim X = 2$. If $\dim X \geq 3$, however, we can
use some approximation arguments, reducing the involved $\opa$-equation to compactly supported forms. As a result we can prove

\newtheorem{extension}[vanish2]{Proposition}
\begin{extension}     \label{extension}  \ \\
Let $\ell \in\N$ be sufficiently large, and assume that $s$ is a $CR$-section of $(N^{1,0}_M)^{\ell}$ of class at least $\mathcal{C}^4$. Then there exists a holomorphic section $\tilde s$ of $\tilde N^{\ell}$ on a tubular neighborhood of $M$ such that $\tilde s_{\mid M} = s$.
\end{extension}

{\it Proof.} First we choose a $\mathcal{C}^4$-extension $s_o$ of $s$ to $\tilde X$ such that $\opa s_o$ vanishes to the third order along $M$, i.e. $\vert\opa s_o \vert_{\omega_o} = O(\vert\rho\vert^3)$. \\

Now we consider an exhaustion of $\tilde X$ by strictly pseudoconvex domains $\Omega_\varepsilon = \lbrace z\in\tilde X\mid \rho^2 (z) > \varepsilon^2 \rbrace$. Moreover, we define the annular domains $D_j = \Omega_{\frac{1}{j}}\setminus \ol \Omega_{\frac{2}{j}}$. \\

Then we choose a sequence of smooth cut-off functions $\chi_j$ with compact support in $\Omega_{\frac{1}{j}}$ such that $\chi_j\equiv 1$ on $\ol \Omega_{\frac{2}{j}}$ and $\vert d\chi_j\vert^2_{\tilde\omega} \leq 1$ (this is possible since $\tilde\omega$ is complete on $\tilde X\setminus M$). Then 
\begin{equation}  \label{f_j}
f_j := \opa(\chi_j\opa s_o)\in L^2_{0,2}(\tilde X\setminus M, \tilde N^{\ell},\tilde\omega)\cap\mathrm{Ker}\opa
\end{equation}
 is compactly supported in $D_j$. \\
Applying Lemma \ref{csupport} yields $u_j\in L^2_{0,1}(\tilde X\setminus M, \tilde N^{\ell},\tilde\omega)$ supported in $D_j$ satisfying $\opa u_j = f_j$ and

\begin{eqnarray*}
 \Vert u_j\Vert^2_{\tilde\omega} & \leq &  C^2 j^2\Vert f_j\Vert^2_{\tilde \omega} \lesssim j^2\Vert\opa \chi_j\wedge\opa s_o\Vert^2_{\tilde\omega}\\
& \leq & j^2\int_{D_j} \vert\opa s_o\vert^2_{\tilde\omega} dV_{\tilde\omega}\leq j^2\int_{D_j}\rho^{-2}
\vert\opa s_o\vert^2_{\omega_o} dV_{\omega_o} \lesssim j^{-2}
\end{eqnarray*}

Now $g_j = \chi_j\opa s_o - u_j$ is $\opa$-closed and supported in $\Omega_{\frac{1}{j}}$, hence compactly supported in $\tilde X\setminus M$. Note that we may view $g_j$ as forms with values in $\tilde N^{\ell}\otimes L^k$ (since $L$ is holomorphically trivial outside a compact $K$ of $\tilde X\setminus M$). By Proposition \ref{vanish2}, there exist $k,N\in\N$ such that we can find solutions $h_j\in L^2_{0,0}(\tilde X\setminus M, \tilde N^{\ell}\otimes L^k, \tilde\omega, -N\varphi)$ satisfying $\opa h_j = g_j$. Hence $g_j\in L^2_{0,1}(\tilde X\setminus M, \tilde N^{\ell}\otimes L^k,\tilde\omega,  -\varphi)\cap\mathrm{Im}\opa$. By Lemma \ref{bottom}, we can therefore find $\tilde h_j\in 
L^2_{0,0}(\tilde X\setminus M, \tilde N^{\ell}\otimes L^k,\tilde\omega, -\varphi)$ satisfying $\opa\tilde h_j = g_j$ and $\Vert \tilde h_j\Vert_{\tilde\omega,-\varphi}\leq C_o \Vert g_j\Vert_{\tilde\omega,-\varphi}$.
But
\begin{eqnarray*}
\Vert g_j\Vert^2_{\tilde\omega,-\varphi} & \lesssim & \Vert \chi_j\opa s_o\Vert^2_{\tilde\omega,-\varphi} + \Vert u_j\Vert^2_{\tilde\omega,-\varphi}\\
& \lesssim & \int_{\Omega_{\frac{1}{j}}} \vert \opa s_o\vert^2_{\tilde\omega} \rho^{-2} dV_{\tilde\omega} + \int_{D_j} \vert u_j\vert^2_{\tilde\omega} \rho^{-2} dV_{\tilde\omega} \\
& \lesssim & \int_{\Omega_{\frac{1}{j}}} \rho^6 \rho^{-4} dV_{\omega_o} + j^2 \Vert u_j\Vert^2_{\tilde\omega} \lesssim 1
\end{eqnarray*}

 Therefore the sequence $(\tilde h_j)$ is bounded in $L^2_{0,0}(\tilde X\setminus M, \tilde N^{\ell}\otimes L^k,\tilde\omega,  -\varphi)$, hence has a subsequence that weakly converges to $ h_o\in L^2_{0,0}(\tilde X\setminus M, \tilde N^{\ell}\otimes L^k,\tilde\omega, -\varphi)$. Since $\opa \tilde h_j = \opa s_o$ on $\Omega_{\frac{2}{j}}$, we must therefore have $\opa h_o = \opa s_o$ in $\tilde X\setminus M$. 
Moreover, since $ h_o\in L^2_{0,0}(\tilde X\setminus M, \tilde N^{\ell}\otimes L^k,\tilde\omega,  -\varphi)$, we have $\int_{\tilde X\setminus M}\vert  h_o\vert^2 \rho^{-2} dV_{\tilde\omega} < +\infty$. This clearly implies that the trivial extension of $ h_o$ to $\tilde X$ satisfies $\opa h_o = \opa s_o$ as distributions on  $\tilde X$ (not only on $\tilde X\setminus M$). Hence $h_o$ is of class at least $\mathcal{C}^4$ by the hypoellipticity of $\opa$, and must therefore vanish on $M$.\\

Thus $\tilde s = s_o -  h_o$ is a holomorphic section of $\tilde N^{\ell}$ in a tubular neighborhood of $M$ (where the line bundle $L$ is holomorphically trivial) extending $s$. \hfill$\square$\\

\newtheorem{csupport}[vanish2]{Lemma}
\begin{csupport}     \label{csupport}  \ \\
Let $\ell \in\N$ be sufficiently large and $f_j$ be defined by (\ref{f_j}).
For some constant $C> 0$, independent of $j\in\N$,  there exists $u_j \in L^2_{0,1}(\tilde X\setminus M, \tilde N^{\ell},\tilde\omega)$, supported in $D_j$, such that $\opa u_j = f_j$ and
$$\Vert u_j\Vert_{\tilde \omega} \leq C j \Vert f_j\Vert_{\tilde \omega}.$$
\end{csupport}

{\it Proof.} 
We will fix a hermitian metric $\tilde h$ on $\tilde N$. By assumption on $N^{1,0}_M$, we may choose $\tilde h$ such that $(\tilde N,\tilde h)$ is positive near $M$.\\

Replacing $\tilde N^\ell$ by $\tilde N^\ell\otimes K_{\tilde X}^\ast =: F = F(\ell)$ (which is still positive for $\ell$ sufficiently large), we may also assume that $f_j$ is an $(n,2)$-form rather than a $(0,2)$-form. \\

Note that the boundary of $D_j$ consists of two parts: a strictly pseudoconvex part $\pa\Omega_{\frac{1}{j}}$ and a strictly pseudoconcave part $-\pa\Omega_{\frac{2}{j}}$. Since $n\geq 3$, this implies that $D_j$ satisfies condition $Z(n-2)$ (see \cite{FK}), hence the $\opa$-Neumann problem satisfies subelliptic estimates in degree $(p,n-2)$ for all $0\leq p\leq n$.  \\

Now we use a duality argument from \cite{Ch-S}: Let $\opa_c$ be the strong minimal realization of $\opa$ on $L^2_{n,1}(D_j, F,\omega_o)$. Then by Theorem 3 of \cite{Ch-S} the range of $\opa_c$ is closed in $L^2_{n,2}(D_j, F, \omega_o)$, and $\opa_c$-exact forms $f\in
L^2_{n,2}(D_j, F, \omega_o)$ are characterized by the usual orthogonality condition:
$$\int_{D_j} f \wedge\eta = 0   \quad\forall \eta \in L^2_{0,n-2}(D_j,F^\ast,\omega_o)\cap\mathrm{Ker}\opa $$
 But, using Stokes' theorem, we get for $\eta \in \mathcal{C}^\infty_{0,n-2}(\ol D_j,F^\ast)\cap\mathrm{Ker}\opa$
$$ \int_{D_j} \opa(\chi_j\opa s_o) \wedge\eta = \int_{\pa D_j}(\chi_j\opa s_o)\wedge\eta = -\int_{\pa\Omega_{\frac{2}{j}}} \opa s_o\wedge\eta = -  \int_{\pa\Omega_{\frac{2}{j}}} \opa (s_o\wedge\eta)  = 0,$$
and this also holds for $\eta \in L^2_{0,n-2}(D_j,F^\ast,\omega_o)\cap\mathrm{Ker}\opa$ using a standard approximation argument and the subelliptic estimates in degree $(0,n-2)$.

 Hence $f_j = \opa(\chi_j\opa s_o)$ belongs to the image of $\opa_c$, i.e. there exists $u_j\in L^2_{n,1}(D_j, F, \omega_o)$ satisfying $\opa_c u_j = f_j$. As usual, we assume that $u_j$ is the minimal $L^2$-solution i.e. $u_j \in L^2_{n,1}(D_j, F,\omega_o)\cap(\mathrm{Ker}\opa_c)^\perp$.
 In particular, $u_j$ is smooth on $\ol D_j$, and the trivial extension of $u_j$ by zero outside $\ol D_j$ (which we still denote by $u_j$), satisfies $\opa u_j = f_j$ as distributions on $\tilde X\setminus M$ (by definition of the strong minimal realization $\opa_c$). It remains to estimate $\Vert u_j\Vert_{\tilde\omega}$.\\

To do so, we may assume that $u_j = \vartheta\alpha_j$ for some $\alpha_j\in L^2_{n,2}(D_j, F, \omega_o)\cap\mathrm{Dom}\vartheta\cap\mathrm{Dom}\opa_c$ satisfying $\opa_c \alpha_j = 0$, i.e. $(\opa_c\vartheta + \vartheta\opa_c)\alpha_j = f_j$. By the subelliptic estimates, $\alpha_j$ is also sufficiently smooth on $\ol D_j$ ($f_j$ is of class $\mathcal{C}^3$ and vanishes outside a compact of $D_j$, so $\alpha_j$ is at least in the Sobolev space $W^3$ and smooth up to the boundary outside the support of $f_j$).\\

We will now estimate $\alpha_j$ by using a priori estimates from \cite{Gri} for negative line bundles over the strictly pseudoconcave domains $W_j = \tilde X\setminus\ol\Omega_{\frac{2}{j}}$. 
From $\cite[\mathrm{Theorem\ VI\ and\ Theorem\ 7.4}]{Gri}$ it follows that there exists $\lambda >0$ such that
$$\Vert v\Vert^2_{\omega_o, W_j} \leq \frac{\lambda}{\ell} (\Vert \opa v\Vert^2_{\omega_o,W_j} + \Vert\opa^\ast v\Vert^2_{\omega_o, W_j})$$
for all $v\in L^2_{0,q}(W_j,  F^\ast,\omega_o) \cap\mathrm{Dom}\opa\cap\mathrm{Dom}\opa^\ast$, $0\leq q\leq n-2$.
From this we infer by Serre duality as in \cite{Ch-S} that
\begin{equation}  \label{estimatewj}
\Vert v\Vert^2_{\omega_o, W_j} \leq \frac{\lambda}{\ell} (\Vert \opa_c v\Vert^2_{\omega_o, W_j} + \Vert \vartheta v\Vert^2_{\omega_o,W_j})
\end{equation}
for all $v\in L^2_{n,q}(W_j, F,\omega_o) \cap\mathrm{Dom}\opa_c\cap\mathrm{Dom}\vartheta$, $q\geq 2$.\\

Note that $D_j$ has two connected components $D_j^\pm$.
We now choose an extension $\tilde\alpha_j^\pm$ of ${\alpha_j}_{\mid D_j^\pm}$ to $W_j$ such that $\tilde\alpha_j^\pm \in \mathrm{Dom}\opa_c\cap\mathrm{Dom}\vartheta$ (on $W_j$!) and such that $$\Vert\vartheta\tilde \alpha_j^\pm\Vert^2_{\omega_o, W_j} + \Vert \opa\tilde \alpha_j^\pm\Vert^2_{\omega_o, W_j} \leq b( \Vert\vartheta \alpha_j\Vert^2_{\omega_o, D_j^\pm} + \Vert \opa \alpha_j\Vert^2_{\omega_o, D_j^\pm} + \Vert\alpha_j\Vert^2_{\omega_o, D_j^\pm})$$ for some constant $b$ not depending on $\alpha_j$ nor on $j$. This is possible for $j$ sufficiently large by general Sobolev extension methods (locally we flatten the boundary $\pa\Omega_{\frac{1}{j}}$ and extend the sufficiently smooth $\alpha_j$ componentwise  across $\pa\Omega_{\frac{1}{j}}$ by first order reflection, then we use a partition of unity, cf e.g. \cite{E}).\\

Applying (\ref{estimatewj}) with $\tilde\alpha_j^\pm$ yields
\begin{eqnarray*}
\Vert \tilde\alpha_j^\pm\Vert^2_{\omega_o, W_j} & \leq & \frac{\lambda}{\ell} (\Vert \opa_c \tilde\alpha_j^\pm \Vert^2_{\omega_o, W_j} + \Vert \vartheta \tilde\alpha_j^\pm\Vert^2_{\omega_o, W_j}) \\
& \leq & \frac{\lambda}{\ell} b ( \Vert\vartheta \alpha_j\Vert^2_{\omega_o, D_j^\pm} + \Vert \opa \alpha_j\Vert^2_{\omega_o, D_j^\pm} + \Vert\alpha_j\Vert^2_{\omega_o, D_j^\pm})\\
& = & \frac{\lambda}{\ell} b ( \Vert\vartheta \alpha_j\Vert^2_{\omega_o, D_j^\pm} +\Vert\alpha_j\Vert^2_{\omega_o, D_j^\pm})
\end{eqnarray*}
For $\ell$ sufficiently large we thus obtain
\begin{eqnarray*}
\Vert\alpha_j\Vert^2_{\omega_o} & \leq &\Vert\vartheta \alpha_j\Vert^2_{\omega_o} = \ \ll \opa_c\vartheta\alpha_j,\alpha_j\gg_{\omega_o} \\
& = & \ll f_j,\alpha_j\gg_{\omega_o} \leq\Vert f_j\Vert_{\omega_o} \Vert\alpha_j\Vert_{\omega_o},
\end{eqnarray*}
which implies
$$\Vert\alpha_j\Vert_{\omega_o} \leq \Vert f_j\Vert_{\omega_o}.$$

Thus
$$\Vert u_j\Vert^2_{\omega_o}\ =\ \ll\vartheta\alpha_j,\vartheta\alpha_j\gg_{\omega_o} \ =\ \ll \opa_c\vartheta\alpha_j, \alpha_j\gg_{\omega_o} \ =\ \ll f_j,\alpha_j\gg_{\omega_o}\ \leq\ \Vert f_j\Vert^2_{\omega_o}.$$
It remains to compare the norms $\Vert u_j\Vert_{\omega_o}$ and $\Vert u_j\Vert_{\tilde\omega}$. To do so, we re-identify $u_j$ and $f_j$ with $\tilde{N}^\ell$-valued $(0,1)$ and $(0,2)$-forms  again. Since $M$ is Levi-flat, we have $dV_{\tilde\omega} \sim \rho^{-2} dV_{\omega_o}$. Using the Levi-flatness of $M$ again, we also have $\vert f_j\vert^2_{\omega_o} = \vert \opa\chi_j\wedge \opa s_o\vert^2_{\omega_o} \lesssim \rho^{-2} \vert f_j\vert^2_{\tilde\omega}$. On the other hand, we have $\tilde\omega \gtrsim \omega_o$, which implies $\vert u_j\vert^2_{\omega_o} \gtrsim \vert u_j\vert^2_{\tilde\omega}$. Since $u_j$ is supported in $D_j$, we thus have
$$\Vert u_j\Vert^2_{\tilde\omega} \lesssim j^2 \Vert u_j\Vert^2_{\omega_o} \leq j^2 \Vert f_j\Vert^2_{\omega_o} \lesssim j^2 \Vert f_j\Vert^2_{\tilde\omega},$$
which proves the desired estimate.
\hfill$\square$\\

The point of the following lemma is that even though $\ell\in\N$ can be arbitrary big, the weight function $-\varphi$ stays the same (it does not have to be multiplied by a large integer as $\ell$ increases!).\\

\newtheorem{bottom}[vanish2]{Lemma}
\begin{bottom}     \label{bottom}  \ \\
Let $\ell,k \in\N$ be arbitrary. Then $\mathrm{Im}\opa$ is closed in $L^2_{0,1}(\tilde X\setminus M, \tilde N^{\ell}\otimes L^k,\tilde\omega, -\varphi)$. This implies that there exists a constant $C_o$ such that every\\ $f\in L^2_{0,1}(\tilde X\setminus M, \tilde N^{\ell}\otimes L^k,\tilde\omega,-\varphi)\cap\mathrm{Im}\opa$ has a solution\\
$u \in L^2_{0,0}(\tilde X\setminus M, \tilde N^{\ell}\otimes L^k,\tilde\omega,-\varphi)$ satisfying $\opa u = f$ and\\ $\Vert u\Vert_{\tilde\omega,-\varphi} \leq C_o \Vert f\Vert_{\tilde\omega,-\varphi}$.
\end{bottom}

{\it Proof.} We will show that for $\ell, k\in\N$ arbitrary, $\mathrm{Im}\opa$ is closed in\\
 $L^2_{n,n}(\tilde X\setminus M, (\tilde N^\ast)^\ell\otimes (L^\ast)^k,\tilde\omega,  \varphi)$; then we argue by duality.\\

First we consider a smooth extension of the the hermitian metric $\omega_o$ on $X$ to $\tilde X$.
Recall that in degree $(n,n)$, the curvature term in (\ref{curvature}) is given by the trace of the curvature form with respect to the metric under consideration. We now modify the metric in $(\tilde N^\ast)\ell\otimes (L^\ast)^k$ by a bounded factor $\exp{(-m\rho^2)}$. This adds to the curvature a term which is $mi\pa\opa\rho^2= 2m\rho i\pa\opa\rho + 2mi\pa\rho\wedge\opa\rho$. Taking $m$ sufficiently large, we may therefore assume that $\mathrm{Trace}_{\omega_o}(i\Theta ((\tilde N^\ast)^\ell\otimes (L^\ast)^k))$ is positive outside a compact of $\tilde X\setminus M$. 
Also, by a theorem of Greene and Wu, $\tilde X\setminus M$ admits a strongly subharmonic exhaustion function with respect to the metric $\omega_o$ (since it is non compact). Pasting a multiple of this exhaustion function together with $\varphi$ (and still calling this modified exhaustion function $\varphi$), we may therefore assume that 
$$\mathrm{Trace}_{\omega_o}(i\Theta((\tilde{N}^\ast)^\ell\otimes (L^\ast)^k) + i\pa\opa\varphi) \gtrsim \rho^{-2}$$
on $\tilde X\setminus M$.
If necessary, we can also modify the metric $\omega_o$ such that its torsion can be absorbed by the right-hand side of the above inequality.
But then the above estimate implies that for $f\in L^2_{n,n}(\tilde X\setminus M, (\tilde N^\ast)^\ell\otimes (L^\ast)^k,\omega_o,  \varphi,\mathrm{loc})$, there exists $u\in L^2_{n,n-1}(\tilde X\setminus M, (\tilde N^\ast)^\ell\otimes (L^\ast)^k,\omega_o, \varphi)$ such that
$$\int_{\tilde X\setminus M} \vert u\vert^2_{\omega_o} e^{-\varphi}dV_{\omega_o} \lesssim 
\int_{\tilde X\setminus M} \rho^2 \vert f\vert^2_{\omega_o} e^{-\varphi}dV_{\omega_o},$$
provided the latter integral is finite.
 Since $u$ is of bidegree $(n,n-1)$ and  $\tilde\omega \gtrsim \omega_o$, we get
$$\int_{\tilde X\setminus M} \vert u\vert^2_{\tilde\omega} e^{-\varphi}dV_{\tilde \omega} \lesssim \int_{\tilde X\setminus M} \vert u\vert^2_{\omega_o} e^{-\varphi}dV_{\omega_o}.$$

On the other hand, the  Levi-flatness of $M$ implies $dV_{\tilde\omega} \lesssim \rho^{-2} dV_{\omega_o}$.
 Since $f$ is of top degree, this means that 
$$\int_{\tilde X\setminus M} \rho^2 \vert f\vert^2_{\omega_o} e^{-\varphi}dV_{\omega_o} \lesssim
\int_{\tilde X\setminus M} \vert f\vert^2_{\tilde\omega} e^{-\varphi}dV_{\tilde \omega}.$$
So we get 
$$\Vert u\Vert_{\tilde\omega,\varphi} \lesssim \Vert f\Vert_{\tilde\omega,\varphi}.$$
This proves that $\mathrm{Im}\opa$ is closed in
 $L^2_{n,n}(\tilde X\setminus M, \tilde N^{\otimes(-\ell)}\otimes L^{\otimes (-k)},\tilde\omega, \varphi)$.\\

Now Serre duality permits to conclude that $\mathrm{Im}\opa$ is closed in\\
 $L^2_{0,1}(\tilde X\setminus M, \tilde N^{\ell}\otimes L^k,\tilde\omega,  -\varphi)$. \hfill$\square$\\

\section{Realization as a hypersurface of a compact K\"ahler manifold}

The final step in the proof of Theorem \ref{main} is to show the following Proposition.

\newtheorem{kähler}{Proposition}[section]
\begin{kähler}     \label{kähler}  \ \\
A tubular neighborhood of $M$ in $X$ can be holomorphically embedded into a compact K\"ahler manifold of dimension $n$.
\end{kähler}

{\it Proof.} In \cite{OS2} (see also \cite{O6}) and \cite{HM}), Kodaira's embedding theorem was generalized to the setting of compact Levi-flat $CR$ manifolds, and it was shown that sufficiently high tensor powers of a positive $CR$ line bundle over a smooth, compact Levi-flat hypersurface $M$ admit enough $CR$ sections $s_0,\ldots, s_m$, so that the $CR$ map $\lbrack s_0 : \ldots : s_m\rbrack$ provides a $CR$-embedding of $M$ into $\C\PP^m$. This applies to our situation, since $N^{1,0}_M = \tilde N_{\mid M}$ is assumed to be positive.\\

In particular, it was proved in $\cite{HM}$ that if $\ell$ is big enough, then the $\mathcal{C}^4$-smooth $CR$-sections of $\tilde N^{\otimes\ell}_{\mid M}$ separate the points on $M$ and give local coordinates on $M$. Using Proposition \ref{extension}, the $CR$-sections of $\tilde N^{\otimes\ell}_{\mid M}$ can be extended to holomorphic sections of $\tilde N^{\otimes\ell}$ over a tubular neighborhood of $M$ in $\tilde X$.\\

Arguing by continuity, it is not difficult to see that if $\ell$ is big enough, then, after possibly shrinking $X$, the holomorphic sections of $\tilde N^{\otimes\ell}$ separate points and give local coordinates on $X$.  Hence we have  a
 holomorphic embedding $\Psi: X \hookrightarrow \C\PP^m$.\\

 Let $\omega_{FS}$ denote the Fubini-Study metric on $\C\PP^m$.  We will show that $\Psi^\ast \omega_{FS}$ extends to a K\"ahler metric on a divisorial blow-up of $\tilde X$. \\

First we extend $\Psi$ to  a meromorphic map $\tilde\Psi: \tilde X\longrightarrow \C\PP^m$: The embedding $\Psi$ is obtained from holomorphic sections $s_j,\ j=0,\ldots,m$ of the line bundle $\tilde N^{\otimes\ell}$ over $X$ for some large $\ell$. Each of these holomorphic sections $s_j$ can be extended to a holomorphic section $\tilde s_j$ of the line bundle $\tilde N^{\otimes\ell} \otimes L^{\otimes k}$ over $\tilde X$ for some large $k$ (use Proposition \ref{vanish2} and see the  proof of Proposition \ref{extension}; also note that $L$ is trivial over $X$).  Then $\lbrack \tilde s_0: \ldots : \tilde s_m\rbrack$ gives a meromorphic extension $\tilde\Psi$ of $\Psi$.\\

By Hironaka's method we may blow up $\tilde X$ along smooth centers, several times if necessary, to obtain a smooth complex manifold $\hat X$ of dimension $n$, together with a holomorphic map $p:\hat X\longrightarrow \tilde X$, such that the induced map $\hat\Psi = \tilde\Psi \circ p : \hat X\longrightarrow \C\PP^m$ is holomorphic. Let $Z$ denote the exceptional divisor of $p$. \\

We  have $\hat\Psi^\ast\omega_{FS} \geq 0$ on $\hat X$, and $\hat\Psi^\ast\omega_{FS} > 0$ on $\lbrace z\in \hat X\mid \mathrm{Jac}\hat\Psi(z)\not= 0\rbrace$. But  since $\Psi$ gives an embedding of $X$,
the analytic set  $\lbrace z\in \hat X\mid \mathrm{Jac}\hat\Psi(z)= 0\rbrace$ is
compact in the strongly pseudoconvex manifold $\hat X\setminus M$. But this means that  $\lbrace z\in \hat X\mid \mathrm{Jac}\hat\Psi(z)= 0\rbrace \subset Z$. \\

We  choose a 
relatively compact open subset $V$ of $\hat X\setminus M$ that contains $Z$.\\

According to Grauert $\cite[\S 3,\ \mathrm{Satz}\ 1]{G1}$, the line bundle $\mathcal{O}(Z)$  associated to the divisor $Z$ is negative.  The curvature form $\Omega = i\Theta (\mathcal{O}(-Z))$ defines a positive K\"ahler form on $V$. Since $\mathcal{O}(-Z)$ is trivial over $V\setminus Z$, there exists a smooth function $\psi$ such that $i\pa\opa\psi = \Omega$ on $V\setminus Z$. Now we choose a smooth cut-off function $\lambda$ with compact support in $V$ such that $\lambda \equiv 1$ in a neigborhood of $Z$. For some sufficiently small $\tau > 0$, the form
\begin{equation} \nonumber
\omega = \left\{ \begin{aligned}
\hat\Psi^\ast\omega_{FS} + \tau i\pa\opa (\lambda\psi)\quad  \mathrm{on}\ \hat X\setminus Z\\
\hat\Psi^\ast\omega_{FS} + \tau \Omega \qquad\qquad\quad \mathrm{near}\ Z
\end{aligned}
\right.
\end{equation}
is then a K\"ahler metric on $\hat X$. \hfill$\square$\\

{\it End of the proof of Theorem \ref{main}.}
The last step in the proof of the theorem is as in \cite{Br} or \cite{O7} (see also section 2). By Proposition \ref{kähler}, $M$ can be realized as a smooth real hypersurface in a compact K\"ahler manifold $\hat X$. Repeating the arguments from section 6, the holomorphic normal bundle $N^{1,0}_M$ (which is topologically trivial) extends to a topologically trivial holomorphic  line bundle over $\hat X$. The Hermitian metric on $N^{1,0}_M$ can be extended to a Hermitian metric of this holomorphic line bundle. 
Since the holomorphic line bundle is topologically trivial, its curvature form is $d$-exact. Applying the $\pa\opa$-lemma on K\"ahler manifolds, we may thus conclude that the curvature form admits a potential, i.e. there exists a smooth function on $M$ which is strictly plurisubharmonic along the leaves of the Levi foliation of $M$. This contradicts the maximum principle. Therefore such $M$ cannot exist. \hfill$\square$\\

\end{document}